\documentclass[12pt]{article}
\author{M.~E.~Kazarian\thanks{Steklov Mathematical Institute
RAS and the Poncelet Laboratory, Independent University of Moscow,
partly supported by the RFBR Grant 04-01-00762}, S.~K.~Lando
\thanks{Institute for System Research RAS and the Poncelet Laboratory,
Independent University of Moscow, partly supported by the grant
ACI-NIM-2004-243 (N\oe ds et tresses)}
\thanks{Both authors are partly supported by the RFBR Grant
05-01-01012-a, NWO-RFBR 047.011.2004.026 (RFBR 05-02-89000-NWOa),
GIMP ANR-05-BLAN-0029-01}}

\title{Thom polynomials for maps of curves with isolated singularities}

\usepackage{amssymb,amsmath}
\usepackage{theorem}

\def\C{{\mathbb C}}

\def\N{{\bf N}}

\def\Q{{\mathbb Q}}

\def\cM{{\cal M}}

\def\cN{{\cal N}}
\def\cO{{\cal O}}
\def\cR{{\cal R}}

\def\deg{{\rm deg}}

\def\codim{{\rm codim}}

\def\Aut{{\rm Aut}}

\def\ve{{\varepsilon}}
\def\pt{{\rm pt}}

\newtheorem{theorem}{Theorem}[section]

\newtheorem{proposition}[theorem]{Proposition}
\newtheorem{corollary}[theorem]{Corollary}

{\theorembodyfont{\rmfamily}    
\newtheorem{remark}[theorem]{Remark}

}

\begin{document}

\maketitle

\begin{flushright}
{\it To Vladimir Igorevich Arnold\\
on the occasion of his 70 birthday\\
with gratitude and admiration}
\end{flushright}

\begin{abstract}
Thom (residual) polynomials in characteristic classes are used in
the analysis of geometry of functional spaces. They serve as a
tool in description of classes Poincar\'e dual to subvarieties of
functions of prescribed types. We give explicit universal
expressions for residual polynomials in spaces of functions on
complex curves having isolated singularities and
multisingularities, in terms of few characteristic classes. These
expressions lead to a partial explicit description of a
stratification of Hurwitz spaces.
\end{abstract}

\section{Introduction}

In~\cite{Ar} V.~I.~Arnold investigated spaces of Laurent (or
trigonometric) polynomials in one variable. Each such space is
determined by a pair of positive integers, the orders of the
poles. Extending earlier results by Looijenga~\cite{Lo} and
Lyashko (see~\cite{Ar1}) valid for polynomials he constructed a
compactification of the space of Laurent polynomials and proved
the $K(\pi,1)$-property for the subspace of nondegenerate Laurent
polynomials (those whose all critical values are pairwise
distinct). The Lyashko--Looijenga map taking each function to the
set of its critical values extends to a polynomial finite map on
Arnold's compactification, and in~\cite{Ar} the degree of this map
was computed.

Degenerate Laurent polynomials form the discriminant, which is
stratified according to the degeneration types. In the present
paper, we deduce formulas for the degree of the Lyashko--Looijenga
map restricted to so-called primitive strata, that is, strata of
Laurent polynomials with a single finite degenerate value.  For
the case of polynomials the degrees of the restriction of the
Lyashko--Looijenga map to all strata were computed
in~\cite{LaZvo99, LaZvo03,Z1,Z2}. In fact, our results are of much
more general nature. They are based on the study of spaces of
meromorphic functions on complex curves (Hurwitz spaces) initiated
in~\cite{KL04}. The main tool of the study is the theory of
universal residual polynomials in characteristic classes developed
mainly by R.~Thom~\cite{Th1} for singularities and
Kazaryan~\cite{Kaz03} for multisingularities. In principle, this
theory allows one to describe (the cohomology classes Poincar\'e
dual to) strata in functional spaces formed by functions with
singularities of prescribed types. However, explicit calculations
of the residual polynomials often prove to be cumbersome.

In~\cite{KL04}, our calculations were mainly based on the
indefinite coefficients method efficiently applied to computing
universal polynomials by Rimanyi~\cite{Rim}. Numerous calculations
of this kind we made led to a number of conjectures concerning the
explicit form of these polynomials. A part of them, concerning the
strata of multisingularities in unfoldings of isolated
singularities, we prove here.

A holomorphic function on a smooth complex curve can have isolated
singularities only of type~$A_n$, that is, those having the form
$z\mapsto z^{n+1}$ in appropriate coordinates. In holomorphic
families whose generic element is smooth, curves with double
points arise in an unavoidable way. A typical example is the
family of hyperbolas $xy=\ve$ on the complex plane, which
degenerates into the pair of coordinate axes at $\ve=0$, each
being a branch of the curve. A function having an isolated
singularity at the double point belongs to the type $I_{k,l}$ if
its restriction to one branch at this point is of type~$A_{k-1}$
and to the other branch is of type~$A_{l-1}$. The types $A_n$ and
$I_{k,l}$ are the only possible types of isolated singularities of
functions on curves, and our goal will be the analysis of their
universal unfoldings.

The universal polynomials for singularities are usually expressed
in terms of the Chern classes of (the tangent bundles over) the
manifolds under study. In particular, the number of basic classes
grows as the complexity, whence the codimension, of the
singularity grows. In the case of families of functions on curves,
however, one can manage with finitely many basic classes, whatever
is the codimension. In particular, if the functions in the family
acquire only isolated singularities, then four basic classes are
sufficient~\cite{KL04}. Specializing the universal polynomials in
the four basic classes to spaces of versal unfoldings of isolated
singularities one can obtain explicit formulas for double Hurwitz
numbers.


\section{Enumeration of Laurent polynomials with
prescribed multisingularities}

\subsection{Stratification of standard versal unfoldings}

Consider the following two families of rational functions:
\begin{equation}\label{A}
 x\mapsto x^{n+1}+a_{2}x^{n-2}+\dots+a_{n}x+a_{n+1}
\end{equation}
\noindent
and
\begin{equation}\label{I}
x\,y=\varepsilon,\quad (x,y)\mapsto
 x^k+y^\ell+a_1x^{k-1}+\dots+a_{k-1}x+c+b_{\ell-1}y+\dots+b_1y^{\ell-1}.
\end{equation}

The domain of a function in the first family is the complex line
with coordinate~$x$ and the coefficients
$(a_2,\dots,a_{n+1})\in\C^{n}$ form the parameter space of the
family. This family of polynomials is the standard versal
deformation of the singularity~$A_n$.

The second family is parameterized by the points
$(\varepsilon,a_1,\dots,a_{k-1},c,b_1,\dots,b_{\ell-1})\in\C^{k+\ell}$.
The domain of a function here is the curve
$\{x\,y=\varepsilon\}\subset\C^2$. This curve is smooth for
$\varepsilon\ne0$ and has a double point for $\varepsilon=0$.
This family is the space of Laurent polynomials of
bidegree $(k,\ell)$, or the standard versal deformation
of the singularity $I_{k,\ell}$.

The zeroes of the differential of a function are called its
critical points. For generic parameter values the functions in
both families have only simple (Morse) critical points with
distinct finite critical values whose number is $n$ and $k+\ell$,
respectively. But for some parameter values the function has more
complicated singularities. We say that a function~$f$ acquires a
local singularity of type~$A_m$ at some point of the source if the
function can be represented in the form $z\mapsto z^{m+1}$ for an
appropriate choice of the local coordinates. A function acquires a
multisingularity of type $A_{m_1,\dots,m_r}$ at a point of the
target if the preimage of this point contains pairwise distinct
points where the function has local singularities of types
$A_{m_1}$, \dots, $A_{m_r}$, respectively. The number of simple
critical values that collapse at a point with multisingularity
$A_{m_1,\dots,m_r}$ is equal to $|m|=m_1+\dots+m_r$. Functions
with multisingularity $A_{m_1,\dots,m_r}$ form the stratum
$\sigma_{m_1,\dots,m_r}$ in the space of functions. The
codimension of this stratum is~$|m|-1$.

For both families, we consider the following problem: \emph{find
the number of generic functions in the family that have critical value $0$
with the multisingularity type $A_{\mu_1,\dots,\mu_r}$, with prescribed
simple critical values}.

These numbers are special cases of what is called \emph{double
Hurwitz numbers}, and a variety of formulas for these numbers
is known~\cite{Ok,GJV}. However, our answer is represented
in a different form and the equivalence of our formulas to the
known ones is by no means evident.

Our approach is close to that of Arnold in~\cite{Ar}. To state it,
let $\mu$ denote the number of parameters in the family. For the
deformation of the singularity~$A_n$ it equals~$n$, and for the
singularity~$I_{k,\ell}$ it is~$k+\ell$. It is equal also to the
number of critical points of a generic function in the family. The
\emph{Lyashko-Loojenga map} $\Lambda:\C^\mu\to \C^\mu$ associates
to a parameter value the unordered tuple of critical values of the
function. This definition is applied if the function has only
simple critical points. However,~$\Lambda$ extends to a
holomorphic (and even polynomial) proper quasihomogeneous map to
the whole space~$\C^\mu$ (see~\cite{Ar,Lo}; explicit coordinate
expressions for the Lyashko--Looijenga map can be found
in~\cite{LZg}, Chapter~5). What we are interested in is actually
the degree of the restriction of $\Lambda$ to a particular
multisingularity stratum. The degree of~$\Lambda$ on the whole
parameter space can easily be computed as the ratio of the
products of coordinate weights in the target and the source, which
yields $(n+1)^{n-1}/(n+1)!$ for the case of polynomials and
$k^k\ell^\ell/(k!\ell!)$ for Laurent polynomials.

In the case of polynomials each multisingularity stratum admits an
explicitly given nonsingular normalization. This allows one to
apply the same quasihomogeneuity argument in order to compute the
corresponding degree.

\begin{theorem}{\bf(\cite{LaZvo99})}
In the case of polynomials, the degree of the restriction of the
Lyashko--Looijenga map to the multisingularity stratum, which is
the closure of functions having multisingularity of type
$m_1,\dots,m_r$, is
$$
(n+1)^{n-1-|m|}\frac{(n-|m|)!}{|\Aut(m_1,\dots,m_r)|(n-r-|m|)!},
$$
where $|\Aut(m_1,\dots,m_r)|(n-r-|m|)!$ denotes the order of the
automorphism group of the tuple $m_1,\dots,m_r$, that is, the
product of the factorials of numbers of coinciding elements in
this tuple.
\end{theorem}

Below, we compute the $A$-contribution to universal formulas using
the tools exploited to obtain this result.

In the case of Laurent polynomials the strata have no good
parametrization. To overcome this difficulty, instead of the
direct geometric study of the stratum we relate the degree of the
restriction of~$\Lambda$ to that stratum to the degree of~$\Lambda$
on the whole space. The most convenient language to
establish a relationship between these degrees is that of \emph{equivariant
cohomology}.

The multiplicative group $\C^*$ of nonzero complex numbers acts of
the spaces of standard versal unfoldings of the
singularities~$A_n$ and $I_{k,\ell}$ by multiplication by an
appropriate power of the complex number. For an element
$\lambda\in\C^*$, this action is
\begin{equation}\label{Aact}
\lambda:(x,a_2,\dots,a_{m+1})\mapsto(\lambda
x,\lambda^2a_2,\dots,\lambda^{n+1}a_{n+1})
\end{equation}
for the unfolding of the singularity~$A_n$ and
\begin{eqnarray}\label{Aact}
&\lambda:(x,y,\ve,a_1,\dots,a_{k-1},c,b_1,\dots,b_{\ell-1})
\nonumber\\&\mapsto (\lambda^\ell x,\lambda^k
y,\lambda^{k+l}\ve,\lambda^\ell
a_1,\dots,\lambda^{(k-1)\ell}a_{k-1},\lambda^{k\ell}c,\lambda^k
b_1,\dots,\lambda^{(\ell-1)k}b_{\ell-1})
\end{eqnarray}
for the unfolding of the singularity $I_{k,\ell}$. The restriction
of this action to the coordinate space $a_i$ (respectively,
$\ve,a_i,c,b_j$) defines an action of $\C^*$ on the parameter
space~$\C^\mu$. Denote by $H_{\C^*}^*(\C^\mu)$ the
$\C^*$-equivariant cohomology of the parameter space~$\C^\mu$ with
respect to this action. Since~$\C^\mu$ is contractible, we have
$H_{\C^*}^*(\C^\mu)\simeq H_{\C^*}^*(\pt)\simeq \Q[\tau]$, where
$\tau=c_1(\cO(1))\in H_{\C^*}^2(\C^\mu)$ is the standard
characteristic class of the $\C^*$-action under consideration
(here and below we work with cohomology with rational
coefficients). We also make use of the similar cohomology rings of
the space of unfolding and its image (both isomorphic
to~$\C^{\mu+1}$. Each of these rings is a ring of polynomials in
one variable~$\tau=c_1(\cO(1))$, and in order to distinguish
between elements of these rings we will sometimes make use of
notation~$\tau_X,\tau_Y$ and~$\tau_B$ for the generators in the
equivariant cohomology rings of the source space (the space of the
unfolding), the target space and the base (the space of the
deformation), respectively.

Consider the stratum of a multisingularity
$\sigma=\sigma_{m_1,\dots,m_r}\subset \C^\mu$. The group~$\C^*$
preserves this stratum, and its projectivization
$P\sigma=(\sigma-\{0\})/\C^*$ is a compact algebraic variety
endowed with a natural fundamental homology class. The number
$$\deg~\sigma=\int_{P\sigma}\tau^{\dim P\sigma}$$ is called the
\emph{degree} of the stratum. Moreover, the equivariant cohomology
class Poincar\'e dual to this subvariety $\sigma\subset\C^\mu$ is
well defined. By definition, this class is proportional
to~$\tau^{|m|-1}$, where $|m|-1=\codim~\sigma=\sum m_i-1$. The
proportionality coefficient can be easily expressed in terms of
the degrees of the stratum and the ambient space:
$$[\sigma]=\frac{\deg~\sigma}{\deg~\C^\mu}\,\tau^{|m|-1}.$$
The degree of~$\sigma$ and of its image, as well as the degree of
the mapping~$\Lambda$ and its restriction $\Lambda|_\sigma$
to~$\sigma$, are related by the following natural equations that
reduce computations of the degree of the map to that of the
degrees of geometric objects:

\begin{proposition}
The degree of the restriction of $\Lambda$ to a stratum $\sigma$,
the degree $\deg~\sigma$ of this stratum and the cohomology class
$[\sigma]$ dual to $[\sigma]$ are subject to the equations
$$\deg~\Lambda|_\sigma=\frac{\deg~\sigma}{\deg~\Lambda(\sigma)}
  =\frac{[\sigma]}{[\Lambda(\sigma)]}
  \frac{\deg~\C^\mu}{\deg~\Lambda(\C^\mu)}
  =\frac{[\sigma]}{[\Lambda(\sigma)]}
  \deg~\Lambda.$$
\end{proposition}

In particular, the generic stratum, which coincides with the
entire deformation space~$\C^\mu$, in the case~$A_n$ has the
degree $1/(n+1)!$, while the degree of its image under~$\Lambda$
is $1/((n+1)^{n}n!)$, and we obtain~$(n+1)^{n-1}$ for the degree
of~$\Lambda$ on this stratum. For the deformation of the
singularity~$I_{k,\ell}$, the degree of the generic stratum is
inverse to the product of the weights of all variables, i.e., it
is  $1/((k+\ell)k!\ell!k^{\ell-1}\ell^{k-1})$, while the degree of
the image is $1/((k+\ell)!(k\ell)^{k+\ell})$, which yields
$$
\deg~\Lambda=(k+\ell-1)!\frac{k^{k+1}\ell^{\ell+1}}{k!\ell!}
$$
for the degree of~$\Lambda$.

\subsection{The degrees of primitive strata in the
standard unfolding of~$I_{k,\ell}$}

Theorem~\ref{thI} in Sec.~3 of the present paper gives universal
formulas for the cohomology classes of primitive strata in the
versal unfolding of the singularity $I_{k,\ell}$. Here we show how
these formulas can be used for computing the degrees of primitive
strata for given orders~$k$ and~$\ell$ of the poles. We restrict
ourselves with the codimension~$1$ strata, the \emph{caustic},
which is formed by functions with multiple critical points, and
the \emph{Maxwell stratum}, which consists of functions having
coinciding critical values at distinct critical points. For strata
of bigger codimension, the computations are similar.

The caustic $\sigma_2\subset B=\C^\mu$ is nothing but the image of
the projection of the stratum $A_2(\C^{\mu+1})\subset
X=\C^{\mu+1}$ to the base. Theorem~\ref{thI} yields for the
$I$-contribution to the dual class the universal expression
$$
A_2(X)=\psi^2-\nu_1\nu_2.
$$
In order to evaluate this expression in the unfolding of the
singularity $I_{k,\ell}$, one must make the substitution
$\psi=k\ell\tau_X,\nu_1=k\tau_X,\nu_2=\ell\tau_X$, whence
$$
A_2(\C^{\mu+1})=((k\ell)^2-k\ell)\tau_X^2.
$$
Note that for $k=\ell=1$, this expression is~$0$, in accordance
with the fact that the caustic in the versal deformation
of~$I_{1,1}$ is empty.

The map~$p$ is the projection to the parameter space. The action
of the Gysin homomorphism $p_*$ on the generator $\tau_X^d$ of the
$2d$-cohomology has the form
$p_*\tau_X^d=\frac{k+\ell}{k\ell}\tau_B^{d-1}$, and we obtain the
following expression for the caustic in the space $B=\C^\mu$:
$$
\sigma_2=(k+\ell)(k\ell-1)\tau_B,
$$
whence the degree of the caustic is the degree of the generic
stratum times $(k+\ell)(k\ell-1)$, i.e.,
$$
\deg~\sigma_2=\frac{k\ell-1}{k!\ell!k^{\ell-1}\ell^{k-1}}.
$$
Dividing this by the degree $1/((k+\ell-2)!(k\ell)^{k+\ell-1})$ of
the image, we obtain as the result the degree of the restriction
of~$\Lambda$ to the caustic:
$$
\deg~\Lambda|_{\sigma_2}=
 (k+\ell-2)!\frac{k^{k}\ell^{\ell}}
 {k!\ell!}(k\ell-1).
$$

For the Maxwell stratum, which is a subvariety in the image of the
unfolding, the universal formula for the $I$-contribution looks
like
$$
A_{1,1}(Y)=\frac{\psi(\nu_1+\nu_2)}{2\nu_1^2\nu_2^2}
(\psi^3(\nu_1+\nu_2)-4\psi^2\nu_1\nu_2+2\nu_1^2\nu_2^2).
$$
Substituting $\psi=k\ell\tau_Y,\nu_1=k\tau_Y,\nu_2=\ell\tau_Y$, we
obtain the expression for this class in the unfolding of the
deformation of~$I_{k,\ell}$,
$$
A_{1,1}(Y)=\frac12k\ell(k+\ell)(k\ell(k+\ell)-4k\ell+2)\tau_Y^2.
$$
This expression vanishes  both for $k=\ell=1$ and for $k=2,\ell=1$
or $k=1,\ell=2$: the Maxwell stratum is empty in the deformations
of the singularities~$I_{1,1}$ and~$I_{1,2}$.

Now, the application of the Gysin homomorphism $q_*$ yields
$$
\sigma_{1,1}=q_*(A_{1,1}(Y))=\frac12(k+\ell)(k\ell(k+\ell)-4k\ell+2)\tau_B,
$$
which gives the following expressions for the degree of the
stratum $\sigma_{1,1}$ and the restriction of~$\Lambda$ to this
stratum:
\begin{align*}
\deg~\sigma_{1,1}&
 =\frac{k\ell(k+\ell)-4k\ell+2}{k!\ell!k^{\ell-1}\ell^{k-1}},\\
\deg~\Lambda|_{\sigma_{1,1}}&=
 (k+\ell-2)!\frac{k^{k}\ell^{\ell}}
 {k!\ell!}\frac12(k\ell(k+\ell)-4k\ell+2).
\end{align*}

\begin{remark}\rm
Note that, similarly to the case of the Maxwell stratum, one could
seek for the expression of the caustic starting with the stratum
$A_2(Y)$ in the target space, the result, of course, would be the
same.
\end{remark}

In the case of a general stratum $\sigma_{m_1,\dots,m_r}$, the
degree of the restriction of $\Lambda$ to this stratum has the
form
$$
\deg\Lambda|_{\sigma_{m_1,\dots,m_r}}=
 (k+\ell-|m|)!\frac{k^{k+2-|m|}\ell^{\ell+2-|m|}}
 {k!\ell!}P_{m_1,\dots,m_r},
$$
where $P_{m_1,\dots,m_r}$ is a symmetric polynomial in~$k$
and~$\ell$. For strata of small codimension, this polynomial is
given by the table below.
$$
{\renewcommand{\arraystretch}{1.5}
\begin{array}{|c|c|}
 \hline (m_1,\dots,m_r)&P_{m_1,\dots,m_r}\\\hline
 (2)&k\ell-1\\
 (1,1)&\frac12(k\ell(k+\ell)-4k\ell+2)\\
 (3)&(k\ell)^2-5k\ell+2(k+\ell)\\
 (2,1)&(k\ell-3)(k\ell(k+\ell)-6k\ell+2(k+\ell))\\
 (1,1,1)&\frac16((k\ell)^2((k+\ell)^2-12(k+\ell)+40)
    +k\ell(6(k+\ell)-80)+24(k+\ell))\\
 \hline
\end{array}}
$$

Note that in order to compute a double Hurwitz number, one has to
divide the degree of~$\Lambda$ restricted to the corresponding
stratum by $k\ell|\Aut(k,\ell)|$. Here the product of~$k$
and~$\ell$ corresponds to the fact that each equivalence class of
a Laurent polynomial has exactly~$k\ell$ representatives in the
standard versal deformation of $I_{k,\ell}$  (see~\cite{Ar}). The
order of the automorphism group is either~$1$ or~$2$ depending on
whether the orders~$k$ and~$\ell$ of the poles coincide or not; if
they coincide, then the transposition of the poles defines an
automorphism of the unfolding.

\section{Universal formulas}
\subsection{Relative characteristic classes and universal polynomials}

The subject of our study are commutative triangles of spaces and
maps of the form
\begin{equation*}
\begin{array}{rcl}
\widetilde X &\stackrel{\tilde f}\longrightarrow&\widetilde Y \\
{}_{\tilde p}\kern-2ex\searrow\kern-2ex&&
\kern-2ex\swarrow\kern-1.5ex{}_{\tilde q}\\
&B
\end{array}
\end{equation*}

Here $\widetilde X,\widetilde Y$, and $B$ are smooth complex
varieties, $\dim~\widetilde X=\dim~\widetilde Y=\dim~B+1$. We
suppose that
\begin{itemize}
\item the fibers of~$\widetilde q$ are smooth complex curves;

\item the fibers of~$\widetilde p$ are nodal compex curves, that
is, their only admissible singularities are points of transversal
double selfintersection (nodes);

\item a section~$\gamma:B\to \widetilde Y$ (fiberwise
``infinity'') is chosen such that the ``poles'' of the
restrictions of~$\tilde f$ to the fibers of~$\widetilde p$ are
smooth points of the fibers, and the orders of these poles are the
same for all fibers;

\item outside of the poles the family of functions is generic (in
particular, if~$\widetilde p$ has singular fibers, then their
images form a subvariety of codimension~$1$ in~$B$).

\end{itemize}

Denote by $Y=\widetilde Y\setminus \gamma(B)$ and $X=\widetilde
X\setminus\widetilde f^{-1}(Y)$ the ``finite parts'' of the
corresponding spaces, and denote by $p,q$ and $f$ the restrictions
of the maps $\tilde p,\tilde q$ and $\tilde f$, respectively, to
the corresponding subspaces. As a result, we obtain the following
diagram:
\begin{equation}\label{eutd}
\begin{array}{rcl}
X &\stackrel{f}\longrightarrow&Y \\
{}_{p}\kern-2ex\searrow\kern-2ex&&
\kern-2ex\swarrow\kern-1.5ex{}_{q}\\
&B
\end{array}
\end{equation}
The maps in this diagram have only generic singularities. The
restriction of~$f$ to a fiber of~$p$ is a holomorphic map of a
nodal curve to a smooth curve, hence the space~$B$ can be treated
as a family of holomorphic maps from nodal to smooth curves,
and~$X$ as the universal curve over~$B$.

The standard unfoldings of the singularity $A_n$ given by
Eq.~(\ref{A}) and $I_{k,\ell}$ given by Eq.~(\ref{I}) provide
examples of such diagrams. In the source of the unfolding $A_n$
all the fibers are smooth; each fiber is the projective line
punctured at infinity. In the source of the unfolding of
$I_{k,\ell}$ the fibers are smooth for $\ve\ne0$, and they are
nodal for $\ve=0$; in the first case, each fiber is the projective
line punctured at two points, while in the second case the fiber
consists of two intersecting projective lines, each punctured at a
point distinct from the point of intersection.

The requirement that the orders of all poles are the same
guarantees that ``the singularities of the fiberwise maps do not
tend to infinity''. This means, for example, that the restriction
of~$p$ to the subvariety~$\Sigma$ of critical points of~$f$ is
proper. Therefore, for any cohomology class in~$X$ supported on
$\Sigma$, its direct image in~$B$ is well defined. To be more
precise, in what follows we must replace everywhere the cohomology
of~$X$ by the relative cohomology $H^*(X,X\setminus \Sigma)$. But
in order to keep notation simple we shall still refer to $H^*(X)$,
having in mind that all the universal polynomials to be considered
are supported on $\Sigma$.

Denote by $\Delta\subset X$ the subvariety of double points
of~$p$; for a generic family~$B$, the codimension of this
subvariety is~$2$. Introduce the relative characteristic classes
of our maps,
\begin{eqnarray*}
\psi&=&f^*(c_1(Y))-p^*(c_1(B))=f^*(c_1(Y)-q^*(c_1(B));\\
\nu&=&c_1(X)-p^*(c_1(B)),
\end{eqnarray*}
and denote by~$\nu_1,\nu_2$ the Chern roots of the normal bundle
to~$\Delta$ in~$X$. The latter are well defined only on~$\Delta$;
they can be thought of as the first Chern classes of the two line
bundles over~$\Delta$ locally given by the tangent lines to the
branches of the curve at the double point. All the formulas for
the universal polynomials depend on the classes $\nu_1,\nu_2$ in a
symmetric way, hence there is no need to ascribe each branch to a
specific class. More precisely,  we use notation
$\nu_1\nu_2\,P(\nu_1+\nu_2,\nu_1\nu_2)$ for cohomology classes
representable in the form $i_*(P(c_1(\N_\Delta),c_2(\N_\Delta)))$,
where $i:\Delta\to X$ is the embedding, and~$P$ is an arbitrary
polynomial in the ordinary Chern classes $c_{1,2}(\N_\Delta)$ of
the normal bundle $\N_\Delta$ to the embedding.

All the four cohomology classes $\psi$, $\nu$, $\nu_1$ and $\nu_2$
can be treated both as elements of the second cohomology group
$H^2(X)=H^2(X,\Q)$, and as elements of the first Chow group
of~$X$. For the sake of definiteness, we speak below only about
the cohomology ring (with rational coefficients) having in mind
that all the results are valid for the Chow ring as well.

The variety~$X$ contains subvarieties $A_m(X)$ consisting of those
smooth points of the fibers of~$p$ where the restriction of~$f$ to
the fiber has a singularity of type~$A_m$. Similarly, the
variety~$Y$ contains subvarieties $A_{m_1,\dots,m_r}(Y)$
consisting of those points whose preimages are smooth points of
the fibers of~$p$, $r$ of them being points where the restriction
of~$f$ to the fiber has singularities of types
$A_{m_1},\dots,A_{m_r}$, while it is nonsingular outside these
points. Denote by $[A_m(X)]\in H^*(X)$, $[A_{m_1,\dots,m_r}(Y)]\in
H^*(Y)$ the cohomology classes Poincar\'e dual to the closures of
these subvarieties. The numbers $m_1,\dots,m_r$ form a partition,
which we shall also write in the form $1^{n_1}2^{n_2}\dots$, where
$n_i$ denotes the number of occurrences of the part~$i$ in the
tuple $m_1,\dots,m_r$; in particular, all but finitely many $n_i$
are~$0$. By definition, $|\Aut(m_1,\dots,m_r)|=n_1!n_2!\dots$.

Applying the general theory of universal polynomials for
singularity classes~\cite{Th1} (respectively, multisingularity
classes~\cite{Kaz03}) to the case under consideration~\cite{KL04}
allows one to assert that the classes $[A_m(X)]$ (respectively,
$[A_{m_1,\dots,m_r}(Y)]$) admit universal expressions in the form
of homogeneous polynomials in the classes $\psi,\nu,\nu_1,\nu_2$
(respectively, in the $f_*$-images of monomials in these classes).
(The set of basic classes we have chosen here differs slightly
from that in~\cite{KL04}. Namely, our choice there was
$\Psi=\psi$, $\Sigma=\psi-\nu$, $\Delta=\nu_1\nu_2$,
$N=-(\nu_1+\nu_2)$.)

\begin{theorem}[\cite{KL04}]\label{texp}
There is a generating function
$$
\cR(t_1,t_2,t_3,\dots)=\sum R_{m_1,\dots,m_r}t_{m_1}t_{m_2}\dots
t_{m_r}
$$
with polynomial coefficients
$R_{m_1,\dots,m_r}=R_{m_1,\dots,m_r}(\psi,\nu,\nu_1,\nu_2)$ such
that for any generic family of holomorphic maps of curves to
curves over a base~~$B$ of the form~{\rm(\ref{eutd})} having only
isolated singularities the classes $[A_m(X)]\in H^*(X)$
{\rm(}respectively, $[A_{m_1,\dots,m_r}(Y)]\in H^*(Y)${\rm)}
coincide with the polynomials $R_m$ {\rm(}respectively, with the
coefficients of $t_1t_2\dots t_r$ in the exponent $\exp f_*(\cR)$
of the direct image of~$\cR${\rm)}.
\end{theorem}

Using the indeterminate coefficients method, we computed
in~\cite{KL04} explicit expressions for universal polynomials of
small codimensions. Now we are going to give formulas for all
universal polynomials, of arbitrary codimension.

\subsection{The formulas}
The restriction of the relative cotangent bundle of~$p$
to~$\Delta$ is trivial. Since all the monomials in~$\nu_1,\nu_2$
arising in the universal formulas are divisible by the
product~$\nu_1\nu_2$ (these monomials make sense only when
restricted to $\Delta=\nu_1\nu_2$), we conclude that the
classes~$\psi,\nu,\nu_1,\nu_2\in H^2(X)$ introduced in the
previous section are subject to the relations
$$
\nu\nu_1=\nu\nu_2=0.
$$
Any monomial
$R\in\Q[\psi,\nu,\nu_1,\nu_2]/\langle\nu\nu_1,\nu\nu_2\rangle$
admits a unique representation in the form
$$
R(\psi,\nu,\nu_1,\nu_2)=R_A(\psi,\nu)+R_I(\psi,\nu_1,\nu_2)-R_0(\psi),
$$
where $R_A(\psi,\nu)=R(\psi,\nu,0,0)$,
$R_0(\psi)=R_A(\psi,0)=R_I(\psi,0,0)$. We call the polynomial
$R_0$ (respectively, $R_A$, $R_I$) the $0$- (respectively, $A$-,
$I$-) {\it contribution\/} to~$R$.

Being elements of the ring
$\Q[\psi,\nu,\nu_1,\nu_2]/\langle\nu\nu_1,\nu\nu_2\rangle$, the
coefficients of the generating function~$\cR$ whose existence is
guaranteed by Theorem ~\ref{texp} split into the $0$-, $A$-, and
$I$-contributions, each of which, in its own turn, is a generating
function:
$$
\cR(\psi,\nu,\nu_1,\nu_2;t)=\cR_A(\psi,\nu;t)
+\cR_I(\psi,\nu_1,\nu_2;t)-\cR_0(\psi;t).
$$
The main result of the paper consists in explicit formulas for
these contributions.

According to the restriction method, in order to compute the
contributions explicitly, it suffices to consider ``elementary
blocks'' in functional spaces, namely, the standard versal
unfoldings of the isolated singularities~(\ref{A}) and~(\ref{I}).

Consider first the standard versal unfolding of the
singularity~$A_n$, that is, in appropriate coordinates
$(x,a_2,a_3,\dots,a_{n+1})$ in~$X$ we have
$$
f:(x,a_2,a_3,\dots,a_{n+1})\mapsto
(x^{n+1}+a_2x^{n-1}+\dots+a_{n+1},a_2,a_3,\dots,a_{n+1})
$$
and both maps~$p$ and~$q$ are the projections to the
$a$-coordinates.
Obviously, our four classes have the following expressions in
terms of the generator~$\tau$ of the ring of equivariant
cohomology: $\nu=\tau$, $\psi=(n+1)\tau$, $\nu_1=\nu_2=0$. Thus,
each polynomial in these classes is reduced to its
$A$-contribution.

\begin{theorem}\label{thA}
For the unfolding of the singularity~$A_n$, the generating
function for the classes of multisingularities is the result of
substitution $\psi=(n+1)\tau$, $\nu=\tau$ to the series
$\cN_A=\cN_A(\psi,\nu;t)$ with rational coefficients given
explicitly by the expansion
$$
\cN_A=1+\psi(\psi-\nu)P_2(\nu,t)+\psi(\psi-\nu)(\psi-2\nu)P_3(\nu,t)+\dots
$$
where the coefficients~$P_m$, $m=2,3,4,\dots$, are given by
\begin{eqnarray*}
1+P_2h^2+P_3h^3+\dots&=&\exp\left(\frac{t_1}\nu h^2+\frac{t_2}\nu
h^3+\frac{t_3}\nu h^4+\dots\right)\\
&=&1+\frac{t_1}\nu h^2+\frac{t_2}\nu
h^3+\left(\frac{t_1^2}{2\nu^2}+\frac{t_3}\nu\right) h^4+\dots
\end{eqnarray*}
\end{theorem}

The terms in the expansion of~$\cN_A$ are graded by either of the
two coinciding gradings. The first grading is given by assigning
degree~$1$ to both basic classes~$\psi$ and~$\nu$. The second one
is obtained by assigning degree~$i$ to the variable~$t_i$ for
$i=1,2,3,\dots$.

Note that the function~$\cN_A$ can be written conveniently as
$$
\cN_A=\exp(t_1\nu d^2/ds^2+t_2\nu d^3/ds^3+t_3\nu d^4/ds^4+\dots)
s^{\frac\psi\nu}|_{s=1},
$$
where~$s$ is an auxiliary variable (cf.~\cite{KL04}).

For the universal unfolding of the singularity~$A_n$ we have
$$
\nu_1=\nu_2=0
$$
and the mapping~$f_*$ acts as follows:
$$
f_*R(\psi,\nu)=(m+1)R((m+1)\tau,\tau)
$$
for any polynomial~$R$, whence

\begin{corollary}
The $A$-contribution~$\cR_A$ to the function~$\cR$ is given by
$$
\frac\psi\nu\cR_A(\psi,\nu;t)=\log~\cN_A(\psi,\nu;t).
$$
\end{corollary}

Although the coefficients of~$P_n$ are rational rather than
polynomial in~$\nu$ (their denominators are powers of~$\nu$), the
coefficients of~$\log~\cN_A$ become polynomial when multiplied
by~$\nu/\psi$, hence the coefficients of~$\cR_A$ are polynomial:
$$
\cR_A=(\psi-\nu)t_1-(\psi-\nu)(2\psi-\nu)t_1^2+(\psi-\nu)(\psi-2\nu)t_2+\dots
$$

Substituting the value $\nu=0$ to the~$A$-contribution, we obtain

\begin{corollary}
The $0$-contribution $\cR_0$ to~$\cR$ is given by
\begin{eqnarray*}
\cR_0(\psi;t)&=&\sum(-2)^{m_1}(-3)^{m_2}\dots(2m_1+3m_2+\dots-1)_{|m|-2}
\psi^{m_1+2m_2+\dots}\frac{t_1^{m_1}}{m_1!}\frac{t_2^{m_2}}{m_2!}\dots\\
&=&\psi\frac{t_1}{1!}+\psi^2\left(\frac{t_2}{1!}-4\frac{t_1^2}{2!}\right)+
\psi^3\left(\frac{t_3}{1!}-6\frac{t_1}{1!}\frac{t_2}{1!}+40\frac{t_1^3}{3!}\right)+\dots,
\end{eqnarray*}
where the summation is carried over all partitions
$1^{m_1}2^{m_2}\dots$. Here $(a)_b$ denotes, for
$b=-1,0,1,2,\dots$, the Pohgammer symbol
$$
(a)_b=a(a-1)(a-2)\dots(a-b+1)=\frac{a!}{b!},\qquad
(a)_{-1}=\frac1a.
$$
\end{corollary}

In order to describe the~$I$-contribution to~$\cR$, let us
introduce another generating function, which we denote by~$\cM_A$,
depending on an auxiliary variable~$z$:
$$
\cM_A(\psi,\nu;z;t)=1+(\psi-\nu)Q_1+(\psi-\nu)(\psi-2\nu)Q_2+
(\psi-\nu)(\psi-2\nu)(\psi-3\nu)Q_3+\dots,
$$
where the rational functions $Q_n=Q_n(\nu;z;t)$ are defined by
means of the expansion
\begin{eqnarray*}
1+Q_1h+Q_2h^2+\dots&=&\exp\left(\frac{t_1}\nu h^2
+\frac{t_2}\nu h^3+\frac{t_3}\nu h^4+\right)/(1-zh)\\
&=&1+zh+\left(\frac{t_1}\nu+z^2\right)h^2+
\left(\frac{t_2}\nu+\frac{t_1}\nu z+z^3\right)h^2+\dots
\end{eqnarray*}

\begin{theorem}\label{thI}
For any $k,\ell\ge1$, the cohomology classes of multisingularities
in the versal unfolding of the singularity~$I_{k,\ell}$ are given
by the generating function obtained as the result of substitution
$\psi=k\ell\tau$, $\nu_1=k\tau$, $\nu_2=\ell\tau$ to the function
$\cN'_I+\cN''_I$, where $\cN'_I$ is given b the formula
$$
\cN'_I(\psi,\nu_1,\nu_2;t)=\cN_A(\psi,\nu_1;t)\cN_A(\psi,\nu_2;t)
$$
and $\cN''_I$ is the result of replacing each monomial of the form
$z^n$ by $(n+1)t_n$ in the product
$$
\psi z\cM_A(\psi,\nu_1;z;t)\cM_A(\psi,\nu_2;z;t).
$$
\end{theorem}

Several first terms of the series~$\cN'_I$ and $\cN''_I$ are
\begin{eqnarray*}
\cN'_I&=&1+\frac{\psi^2(\nu_1+\nu_2)-2\psi\nu_1\nu_2}{\nu_1\nu_2}t_1\\
&&+\frac{\psi^3(\nu_1+\nu_2)-6\psi^2\nu_1\nu_2+2\psi\nu_1\nu_2(\nu_1+\nu_2)}{\nu_1\nu_2}t_2\\
&&+\frac{\psi^4(\nu_1+\nu_2)^2-8\psi^3\nu_1\nu_2(\nu_1+\nu_2)
+24\psi^2\nu_1^2\nu_2^2-6\psi\nu_1^2\nu_2^2(\nu_1+\nu_2)}{2\nu_1^2\nu_2^2}t_1^2+\dots
\end{eqnarray*}
and
\begin{eqnarray*}
\cN''_I&=&2\psi t_1+(6\psi^2-3\psi(\nu_1+\nu_2))t_2\\
&&+\frac{2\psi^3(\nu_1+\nu_2)-12\psi^2\nu_1\nu_2
+4\psi\nu_1\nu_2(\nu_1+\nu_2)}{\nu_1\nu_2}t_1^2+\dots
\end{eqnarray*}

Similarly to the function~$\cN_A$, the terms of the expansions
of~$\cN'_I$ and~$\cN''_I$ are graded by one of the two coinciding
gradings. The first of them is obtained by assigning degree~$1$ to
the classes~$\psi$ and~$\nu_{1,2}$. The second one is the result
of assigning degree~$i$ to the variables~$t_i$, for
$i=1,2,3,\dots$.

Making use of the equation $\nu=0$ valid for the standard
unfolding of $I_{k,\ell}$ and the fact that the action of~$f_*$
has the form
$$
f_*R(\psi,\nu_1,\nu_2)=(k+\ell)R((k+\ell)\tau,k\tau,\ell\tau)
$$
for any polynomial~$R$, we obtain

\begin{corollary}
The $I$-contribution~$\cR_I$ to the function~$\cR$ is given by
$$
\frac{\psi(\nu_1+\nu_2)}{\nu_1\nu_2}\cR_I= \log(\cN'_I+\cN''_I).
$$
\end{corollary}

The coefficients of the function~$\cR_I$ are polynomial rather
than rational as well:
\begin{eqnarray*}
\cR_I&=&1+\psi t_1+\psi^2(t_2-2t_1^2)+\dots\\
&&+\nu_1\nu_2\left(t_1^2-t_2+\frac43(10(\nu_1+\nu_2)-3\psi)t_1^3
-6(\nu_1+\nu_2-3\psi)t_1t_2\right.\\
&&\left.+(2(\nu_1+\nu_2)-5\psi)t_3+\dots\phantom{\frac43}\right);
\end{eqnarray*}
here we separated the $0$-contribution and the terms supported
on~$\Delta$.

\medskip
For the case of monosingularities, the general formulas look like
follows.

\begin{corollary}
The generating function
$$
\sum_{i=0}^\infty[A_i(X)]z^i
$$
for classes of monosingularities in families with isolated
singularities has the form
\begin{multline*}
\sum_{i=0}^\infty[A_i(X)]z^i
 =L(\psi,\nu;z)+\\
 \frac{\nu_1\nu_2}{\nu_1+\nu_2}
\left(\frac{L(\psi,\nu_1;z)}{\nu_1}+\frac{L(\psi,\nu_2;z)}{\nu_2}
+(z^2L(\psi,\nu_1;z)L(\psi,\nu_2;z))'\right)-\frac1{1-\psi z},
\end{multline*}
where
$$
L(\psi,\nu;z)=1+(\psi-\nu)z+(\psi-\nu)(\psi-2\nu)z^2+
(\psi-\nu)(\psi-2\nu)(\psi-3\nu)z^3+\dots
$$
is the $A$-contribution, and $\frac1{1-\psi z}$ is the
$0$-contribution to the corresponding classes.
\end{corollary}

\section{Proofs}
\subsection{Proof of Theorem~\ref{thA}}

The proof of Theorem~\ref{thA} below follows, essentially, the
argument in~\cite{LaZvo99}, see also~\cite{Kaz02}. In order to
compute the cohomology class in
$H^*(Y)=H^*_{\C^*}(\C^{n+1})=\Q[\tau]$ dual to the stratum
$A_{m_1,\dots,m_r}(Y)$, let us construct an explicit
parametrization of this stratum. The points in the closure of
$A_{m_1,\dots,m_r}(Y)$ admit a representation as the value and the
coefficients of the polynomial
$$
P(x)=\frac1{n+1}\int_0^x(\xi-x_1)^{m_1}\cdot\dots\cdot(\xi-x_r)^{m_r}
Q(\xi)d\xi.
$$
Here~$Q$ is a polynomial of an appropriate degree with leading
coefficient~$1$, whose second coefficient is chosen so as to make
the second coefficient in the expanded integrand vanish.

The explicit expression for the integral determines a
 $\C^*$-equivariant map $\rho:\C^{n+1-|m|}\to\C^n$, where
$|m|=m_1+\dots+m_r$ is the codimension of the stratum under
consideration in~$Y$, and the points $x_1,\dots,x_r$ and the
coefficients of~$Q$ serve as coordinates in the domain. This
mapping is finite (in particular, it is proper), its image is
exactly the closure of $A_{m_1,\dots,m_r}(Y)$, and the degree of
the mapping is $|\Aut(m_1,\dots,m_r)|$, the product of the
factorials of the numbers of coinciding parts. Therefore,
$$
\rho_*(1)=|\Aut(m_1,\dots,m_r)|[A_{m_1,\dots,m_r}(Y)].
$$

In order to compute~$\rho_*$, note first that $\rho(0)=0$.
Considering the Poincar\'e dual classes, we see that
$\rho(e(\C^{n+1-|m|}))=e(\C^n)$. Here $e(\cdot)$ denotes the Euler
class of the normal bundle to the origin in the corresponding
space. This class is equal to the product of the weights of all
coordinates. Now the equations $\psi=(n+1)\tau$, $\nu=\tau$, yield
\begin{eqnarray*}
\rho_*(1)
&=&\frac{e(\C^n)}{e(\C^{n+1-|m|})}=(n+1)n(n-1)\dots(n-|m|+2)\tau^{|m|}\\
&=&\frac{\psi(\psi-\nu)\dots(\psi-|m|~\nu)}{\nu^{|m|}}.
\end{eqnarray*}

Taking these classes for the coefficients of the generating
series, we obtain the function~$\cN_A$ in Theorem~\ref{thA}. The
theorem is proved.

\vspace{.5cm}

This argument can be easily modified to construct the generating
function~$\cM_A$ for cohomology classes of multisingularities with
a distinguished singular point. We shall need this generating
function below in the proof of Theorem~\ref{thI}.

Consider the subvarieties $A_{i;m_1,\dots,m_r}(X)$ in the source
space $\C^{n+1}$ consisting of points $(x_0,a_2,\dots,a_{n+1})$
such that the polynomial
$$
x^{n+1}+a_2x^{n-1}+\dots+a_{n+1}
$$
has a singularity of type~$A_i$ at~$x_0$, (pairwise distinct)
critical points of types~$A_{m_1},\dots,A_{m_r}$, and its value at
each of these points coincide with that at~$x_0$. Set
\begin{eqnarray*}
\cM_A(\psi,\nu;z;t_1,t_2,\dots)&=&\sum_{i,m_1,m_2,\dots}
[A_{i;m_1,\dots,m_r}(X)]z^it_{m_1}\dots t_{m_n}.
\end{eqnarray*}

The parametrization similar to the one above proves that the
function~$\cM_A$ coincides with the function in Theorem~\ref{thI}.

\subsection{Proof of Theorem~\ref{thI}}

Now we want to compute the cohomology class
$[A_{m_1,\dots,m_r}(Y)]$ in $\C^*$-equivariant cohomology of the
versal unfolding of the singularity $I_{k,\ell}$. Denote by~$\tau$
the generator of the cohomology ring $H^*(Y)$.

In contrast to the case of~$A_n$ singularity, now the subvariety
$A_{m_1,\dots,m_r}(Y)$ admits no parametrization, and we must
proceed differently. Instead of computing the class dual to the
closure of the subvariety $A_{m_1,\dots,m_r}(Y)$, we compute, in
two different ways, the intersection of this class with the
hypersurface $\ve=0$ consisting of the images of the singular
fibers of~$p$. Comparing the results of the computation, we shall
obtain the desired conclusion.

On one hand, we have
$$
[\ve=0][A_{m_1,\dots,m_r}(Y)]=(k+\ell)\tau[A_{m_1,\dots,m_r}(Y)]
$$

For the purpose of the second computation, let us split the
intersection $A_{m_1,\dots,m_r}(Y)\cap\{\ve=0\}$ into irreducible
components and compute their multiplicities.

As projective line degenerates acquiring a double point, the
monosingularities on this line are distributed somehow between the
two irreducible components of the singular curve. One of the
monosingularities can find itself at the origin. The summand
$\cN'_I$ in the $I$-contribution describes the situation where
there is no singularity at the origin, with the distinguished
critical value, while the summand~$\cN''_I$ is in charge of such a
singularity.

If the function has no singularity at the origin, then its
monosingularities can be split between the two components in all
possible ways, and the restriction to each component is described
by the generating function~$\cN_A$. Each of such splittings forms
an irreducible component of the intersection, and, obviously, the
multiplicity of this component is~$1$. Replacing the class $\nu$
in two copies of the function~$\cN_A$ by the classes $\nu_1$ and
$\nu_2$, respectively, and multiplying the results, we obtain the
contribution of such splittings.

And if there is a monosingularity at the origin, then its
restriction to each of the branches produces two local
singularities. This splitting type is described by the generating
function~$\cM_A$, whose argument~$z$ is in charge of the
distinguished critical point, the one at the origin. Two copies of
this function give the product
$$
\cM_A(\psi,\nu_1;z;t)\cM_A(\psi,\nu_2;z;t).
$$

A function having local singularity~$x^a$ on one branch of the
curve at the origin and $y^b$ on the other branch, can result from
a degeneration of an~$A_m$ singularity on the smoothened curve iff
$a+b=m+1$. Any such partition $(a,b)$ of~$m+1$ determines an
irreducible component of the local singularity. The multiplicity
of intersection of such a component with the divisor $\{\ve=0\}$
is independent of $a$ and $b$ and equals $m+1$. Indeed, this
multiplicity coincides with the dimension of the local algebra
$\C[x,y]/\langle x^a+y^b,xy\rangle$, which is $a+b=m+1$.
Therefore, the multiplicities of all the components are the same,
and in order to compute their contribution one must multiply the
product of the generating functions~$\cM_A$ by~$\psi z$ and
replace in the result each monomial $z^{m+1}$ by $(m+1)t_{m+1}$.
This completes the proof of the theorem.

\end{document}